\documentclass[10pt,a4paper]{amsart}
\input amssym.def
\input amssym.tex
\title[ Hecke groups ]
{ The commutator subgroup and the index formula of the Hecke group  $H_5$}

\author{\SMALL Cheng Lien Lang}
\author{ \SMALL  mong lung lang}

\begin{document}

\baselineskip=12pt


\maketitle

\vspace{-0.3in}

\begin{abstract}
Let
  $A$ be an ideal of
 $\Bbb Z[2cos(\pi/5)]$ and let $  H(A)$ be the principal congruence subgroup of level
  $A$ of the Hecke group $H_5$.
 The present article
  gives explicit formula of
  $[H_5 : H(A)]$. As a byproduct, we prove that the commutator subgroup of
   $H_5$ is not congruence. Consequently, the fact that the commutator subgroup of
    the modular group is congruence cannot be carried over to the Hecke groups.

\end{abstract}

\section{Introduction}
\subsection{} Hecke groups.
 Let $q\ge 3$ be a a fixed integer. The (homogeneous) Hecke group
 $H_q$ is defined to be the maximal discrete subgroup
of  $SL (2,\Bbb R)$
 generated by  $S$ and $T$, where $\lambda _q =2$cos$\,(\pi/q)$,
 $$ S =
\left (
\begin{array}{rr}
0 & 1 \\
-1 & 0 \\
\end{array}
\right ) \,,\,\,
  T = \left (
\begin{array}{rc}
1 & \lambda  _ q\\
0 & 1 \\
\end{array}
\right ) \,.\eqno(1.1)
$$

\noindent

\noindent
Let $A$ be an ideal
 of $\Bbb Z[\lambda_q]$.
We define three  congruence subgroups of $H_q$ as follows.
 $$H_0(A)=\{ (a_{ij})\in H_q \,:\, a_{21} \in A\},\eqno(1.2)$$
 $$
 H_1(A)=\{ (a_{ij})\in H_q \,:\,a_{11}-1, a_{22} -1, a_{21} \in A\},\eqno(1.3)$$
 $$H( A)=\{ (a_{ij})\in H_q \,:\,a_{11} - 1, a_{22} -1, a_{12}, a_{21} \in A\}.\eqno(1.4)$$

 \noindent
 Let $Z_2 = \left < \pm I\right >$. The (inhomogeneous) Hecke group and its congruence subgroups
  are defined as $G_q = H_q/Z_2,$  $ G_0(A) = H_0(A)/Z_2,$ $ G_1(A) = H_1(A)Z_2/Z_2$,
   and $ G(A) = H(A)Z_2/Z_2$.
The present article determines the index $[H_5 : H(A)]$.

\subsection{} The obstruction.
$H(A)$ is the kernel of the  homomorphism $ f: H_5 \to SL(2, \Bbb Z[\lambda]/A)$.
 In the case $A= (2)$, it is known that
 $ f: H_5 \to SL(2, \Bbb Z[\lambda]/2) \cong A_5$ and that
 $[H_5 : H(2)] = 10$ (see [LLT2]). Consequently, the map $f$ is not always
  surjective and
  the finding of the index $[H_q : H(\pi)]$ is
    more involved.
    As a matter of fact, one will see in Section 6 that $f$ is surjective if and only if
 the norm of $A$ is prime to 6. Note that our result is in line with the
  modular group case only if gcd$\,(N(A), 6 ) =1$ (see (1.5)).

\subsection{} The plan and the main results. We consider the
decomposition   $A=\prod A_i$ where $N(A_i)$ is a power of
some rational prime $p_i$ and
gcd$\,(N(A_i), N(A_j)) =1$. To determine the index, we will prove
 in the following sections that
\begin{enumerate}
\item[(i)]$[H_5 : H(A)]=\prod [H_5:H(A_i)]$ (Section 4),
\item[(ii)] $[H_5 :H(A_i)]$ can be expressed in terms of the
 norm $N(A_i)$ (Section 5).
\end{enumerate}
 (i) of th above is achieved by studying the normal subgroup
 $H(A_i)H(A_j)$. (ii) can be done by studying the group structure
of  $SL(2, \Bbb Z[\lambda]/A_i)$.
Our main result shows that   if gcd$\,(N(\pi), 6)=1$, then
   $$[H_5 : H(2^a3^b\pi)]= I_aJ_b\, N(\pi)^3
\prod_{P|(\pi)} (1-N(P)^{-2},\eqno(1.5)$$
where
the product is over the set of all prime ideals that divide $\pi$,
 $I_0 = 1,$ $I_1= 10$, $I_a = 5\cdot 2^{6(a-1)}$ if $a \ge 2$, $J_0=1$,
 $J_b = 120\cdot 3^{6(b-1)}$ if $b \ge 1$.
Our proof is elementary but involves some curial construction of matrices
 (see subsections 3.3 and 3.4). Note that our strategy works
  for  other  groups as well (see [LL2] for elementary matrix groups of P. M. Cohn [C]).

\subsection{} First  application. Subgroups of $H_q$ are called
 {\em congruence} if they contain some principal congruence subgroups.
It is well known that the commutator subgroup of the modular group $SL(2, \Bbb Z)$
 is congruence of level 12. As an application of the index formula given in (1.5),
  we  prove that (unlike the modular group case)  the commutator subgroup of $H_5$ is not congruence (Section 7).
  Two subgroups $A$ and $B$ of $SL(2, \Bbb R)$ are {\em commensurable} with each other if $A\cap B$ is of
   finite index in both $A$ and $B$.
    As   the commutator subgroups of $SL(2, \Bbb Z)=H_3, H_4$ and $H_6$
    are congruence (see Appendix C), our study suggests the following conjecture.

  \smallskip
  \noindent
   {\bf Conjecture.} The commutator subgroup of $H_q$ is congruence if and
 only if $H_q$ is commensurable with the modular group.

\smallskip
A weaker version of the above conjecture is

 \smallskip
  \noindent
   {\bf Conjecture.} Let $q\ge 5$ be a prime.
   Then the commutator subgroup of $H_q$ is not congruence.

\subsection{} Second application. Let $f$ be a function (form) defined on $\mathcal H/\Gamma (n)$,
 where $\mathcal H $ is the union of the upper half plane and $\Bbb Q \cup \{\infty\}$ and
  $\Gamma (n)$ is the principal congruence subgroup of the modular group $\Gamma
   = PSL(2, \Bbb Z)$. Then $\sum f|_{\sigma}$
   is a function (form) of $\mathcal H/\Gamma$,
    where the sum is over the set
    of a complete set of coset representatives of $\Gamma /\Gamma(n))$.
     With the help of our proof of the index formula (1.5), a complete
      set of coset representatives of $G_5/G(\pi)$ can be determined.
       As a consequence, the above construction for functions (forms)
        of $\mathcal H/\Gamma$ can be extended to $G_5$ as well. Note that the space
         of forms for $G_q$ has positive dimension (Theorem 3 of Ogg [O]).

\subsection{}
The remaining of this article is organised as follows. In Section 2, we give
 results related to $SL(2,\Bbb Z[\lambda_q])$ and $[H_5: H(A)]$
where $A$ is a prime ideal. Section 3 is devoted to
  the study of some technical lemmas.  In particular, we give a set of
  generators  for $G(2)$ (Lemma 3.1) which will be used in the
   determination of the index $[H(2) :H(4)]$.  Section 4
    investigates
    the index formula for $H(A)$.
    It is shown that the index formula for the principal congruence subgroup (1.4)
     is multiplicative.
    Section 5  gives the  closed form of the index
       $[H_5: H(A_i)]$.
  The main result of this article can be found in
 Section 6. Throughout the sections, we set $\lambda = 2cos\, \pi/5$.

\section {$SL(2,\Bbb Z[\lambda_q])$}
\subsection{}For $A$ an ideal of $\Bbb Z[\lambda_q]$,
we may also define the  principal congruence subgroup $L (A)$
 of $L_q = SL(2,\Bbb Z[\lambda_q])$ analogously.
The formula for the index of the principal  congruence subgroup in $L_q$
 is easily calculated as the modular group case (see [Sh]);
$$[L_q : L(A)] =
N(A)^3\prod_{P|A}(1-N(P)^{-2}),\eqno(2.1)$$

\noindent where $ N(A)$ denotes the absolute norm of  $A$ in $\Bbb Z[\lambda_q]$
and the product is over the set of all  prime ideals $P$ that divide $A$.
 Let $\pi\in \Bbb Z[\lambda_q]$.
It follows easily from (2.1) that
$$[H(\pi^n) : H(\pi^{n+1})] \le N(\pi)^3.\eqno (2.2)$$

\subsection{} Recall that $\lambda = 2cos\,\pi/5$.  Let $\pi$ be a prime ideal of $\Bbb Z[\lambda]$. Denoted by $p$ the smallest
positive rational prime in $\pi$. Applying our results in
[LLT2], one has
(i) $[H_5 : H(\pi)]= 10$ if $p=2$, (ii)
$[H_5 : H(\pi)]=120$ if $p=3$ or 5
(iii) $[H_5: H(\pi)]= (p-1)p(p+1)$ if $p = 10k \pm 1$,
(iv)
$[H_5: H(\pi)]=(p^2-1)p^2(p^2+1)$ if $3\ne p = 10k\pm 3$.

\section{Technical lemmas }
\subsection{Results of Kulkarni.}In [K], Kulkarni applied a combination of geometric
 and arithmetic methods to show that one can produce
 a set of independent generators in the sense of
 Rademacher for the congruence subgroups of the modular group,
 in fact for all subgroups of finite indices.
 His method can be generalised to all subgroups of
 finite indices of the Hecke groups $G_q = H_q/ Z_2$, where $q$ is a prime. See [LLT1] for detail
  (Propositions 8-10 and section 3 of [LLT1]). By Section 3 of [LLT1], we have the following
   lemma.

\smallskip
\noindent {\bf Lemma 3.1.} {\em The principal congruence subgroup $G(2)=H(2) Z_2/Z_2$ has index $10$ in
 $G_5 H_5/Z_2$. A set of  generators is given by
{\small $$\Omega_2 =\left \{\left (\begin{array}{cc}
   1 & 2\lambda\\
   0 & 1
 \end{array}\right ),
 \left (\begin{array}{cc}
   1 & 0\\
   2\lambda & 1
 \end{array}\right ),
 \left (\begin{array}{cc}
   1+2\lambda & 2+ 2\lambda\\
   2\lambda & 1 +2\lambda
 \end{array}\right ),
 \left (\begin{array}{cc}
   1+2\lambda & 2\lambda\\
   2 +2\lambda & 1+2\lambda
 \end{array}\right )\right \}.\eqno(3.1)$$}}

\subsection{Reduced form.} $a/b$ is in reduced form  if
 {\tiny $\left (\begin{array}{c}
   a\\
   b
 \end{array}\right )$} is a column vector of some $\sigma\in H_5$.
  For any $a, b\in \Bbb Z[\lambda]$ such that the greatest common divisor
   of $a$ and $b$ is a unit, applying results of Leutbecher ([L1, L2]),
    there exists a unique  $n \in \Bbb Z$ such that $a\lambda^n/b\lambda^n$ is in reduced
     form.
    we shall now give an algorithm that enables us
to determine the  reduced factor $n = e(a/b)$.
\smallskip

\noindent
Let $a, b \in \Bbb Z[\lambda] \setminus \{ 0\}$
be given such that
the greatest common divisor of $a$ and $b$  is a unit.
 Then there exists
 a unique rational integer $q$ such that

\begin{enumerate}
\item[(i)] $a = (q \lambda )b + r$,

\item[(ii)]$ -|b\lambda|/2 < r \le |b\lambda|/2$.
\end{enumerate}
We call such a division algorithm {\bf pseudo-Euclidean} (see [R] for more details). In terms of
 matrices, the above can be written as
 $$
\left (
\begin{array}{rr}
1 & -q\lambda \\
0 & 1 \\
\end{array}
\right )
\left (
\begin{array}{r}
a  \\
b \\
\end{array}
\right )  =
\left (
\begin{array}{r}
r  \\
b \\
\end{array}
\right )\,.\eqno(3.2)$$
Note that
$\left (
\begin{array}{rr}
1 & -q\lambda \\
0 & 1 \\
\end{array}
\right ) \in H_5$.
 Applying the  pseudo
 Euclidean algorithm repeatedly, one has,
$$a = (q_1 \lambda) b + r_1\,,$$
$$b = (q_2 \lambda) r_1 + r_2\,,$$
$$r_1 = (q_2 \lambda) r_2 + r_3\,,$$
$$..........................$$
$$..........................$$
$$r_{k+1} = (q_{k+2} \lambda) r_{k+2} + 0\,.$$
The finiteness of the algorithm is
 governed by the fact that
the set of cusps of $H_5$ is $\Bbb Q\,[\lambda]
 \cup \{ \infty \}$.
Note that in terms of matrices, the above can be written as
$$
\left (
\begin{array}{c}
a  \\
b \\
\end{array}
\right )  = A
\left (
\begin{array}{c}
r_{k+2}  \\
0 \\
\end{array}
\right )\,, \eqno (3.3)$$
 where $A \in H_5$.
It is clear that
$\mbox{ gcd}(a,b) =\mbox{ gcd}(b, r_1)
 = \cdots = \mbox{ gcd}
(r_{k+1}, r_{k+2}) =  r_{k+2}$
 is a unit ($a$ and $b$ are coprime).
As  $\lambda$ is a primitive unit,
 there exists $e(a/b) \in \Bbb  Z\cup \{0\}$ such that
 $r_{k+2} = \lambda^{-e(a/b)}\,.$
Multiplying (3.3)  by $\lambda ^{e(a/b)}$,
 one has
$$
\left (
\begin{array}{c}
a\lambda ^{e(a/b))}  \\
b\lambda ^{e(a/b)} \\
\end{array}
\right )  = A
\left (
\begin{array}{c}
\pm1  \\
0 \\
\end{array}
\right )\,.\eqno(3.4)$$
Since $A \in H_5$
 and $\pm 1/0$ is a reduced form,
 we conclude that $
a\lambda ^{e(a/b)}/b\lambda ^{e(a/b)} $ is the reduced form of
 $a/b$.
 We state without proof the following lemma (Proposition 6 of [LLT1]).

 \smallskip
 \noindent {\bf Lemma  3.2.} {\em
{\tiny $\left (\begin{array}{cc}
   a& b\\
   c&d
 \end{array}\right ) $} $\in H_5$ if  and only if
 $ad-bc=1$and
  $a/c$, $b/d$ are in reduced
  forms.}

   \smallskip

 \noindent {\bf Lemma  3.3.} {\em
 Let $p \in \Bbb N$ be an odd prime and let $a = 2\lambda^2$, $c= p\lambda^3$.
 For every  $m \in \Bbb Z$,
    $$  S = \left (\begin{array}{cc}
   1-acm\lambda & a^2m\lambda\\
   -c^2m\lambda & 1+acm\lambda
 \end{array}\right ) \subseteq H(m).\eqno(3.5)$$}

 \smallskip
 \noindent {\em Proof.}  The reduced form of $2/p\lambda $ is
  $2\lambda^2/ p\lambda ^3$. By the definition of reduced form,
   there exists some $u, v \in \Bbb Z[\lambda]$ such that
  $$ X = \left (\begin{array}{cc}
  2\lambda^2  & u\\
  p\lambda ^3 & v
 \end{array}\right )\in H_5.\eqno(3.6)$$
 As a consequence, $ X T^m X^{-1} \in H_5$ for all $m \in \Bbb Z$.
  Note that the matrix form of  $ X T^m X^{-1}$  is given by (3.5).\qed

 \subsection{Some matrices of $H_5$}

 By Lemma 3.2, one can show easily that $H_5$ contains the following matrices
  which will be used for our study.
 {\small $$ \left (\begin{array}{cc}
   1+2\lambda & 2+ 2\lambda\\
   2\lambda & 1+2\lambda
 \end{array}\right ),
 \left (\begin{array}{cc}
   \lambda &\lambda +2\\
   \lambda & 2\lambda +1
 \end{array}\right ),
  \left (\begin{array}{cc}
   -\lambda & \lambda\\
   -2\lambda & 2\lambda +1
 \end{array}\right ),
  \left (\begin{array}{cc}
   3\lambda+2 & -2\lambda-3\\
   4\lambda+3 & -4\lambda-2
 \end{array}\right )\in H_5.\eqno(3.7)$$}

\subsection{}A key proposition (Proposition 4.3) in our  study of the index formula
requires the existence of $\sigma \in H(p)$ such that
  $\sigma \equiv $  {\tiny $
 \left (\begin{array}{cl}
                               1 &   p   \\
                               0 & 1\\
                             \end{array}\right )$} $(mod\, p^{2})$.
                              The matrices we introduced in 3.3 will be used in
                               Section 5 for that purpose. It is perhaps worthwhile
                                to point out that we obtain the above matrices not completely random but by
                                 the study of the fundamental domains of the subgroups
                                  $G(p)$. See [LLT1] for more detail.

\section{$[H_5 : H( \tau)]$ is multiplicative }
In this section,  we  investigate the index formula for the principal
 congruence subgroup (1.4)  and prove that the index formula is multiplicative (Lemma  4.1).
 Note that the index formula for the inhomogeneous Hecke group $G_q$ is not multiplicative.

\smallskip
\noindent {\bf Lemma 4.1.}
{\em  Let  $\tau$, $\pi \in \Bbb Z[\lambda]$ and let $a $ and $b$
 be the smallest positive rational integers in  the ideals $(\tau)$ and
  $(\pi)$ respectively. Suppose that gcd$\,(a, b) =1$.
   Then $H_5= H(\tau)H(\pi)$
   and    $[H_5 : H(\tau  \pi)]
 = [H_5 :H(\tau)][H_5 : H( \pi)]$.}

 \medskip

 \noindent {\em Proof.} Let $K = H(\tau) H(\pi)$.
It is clear that $$     \left (\begin{array}{cc}
                               1 &    a\lambda      \\
                               0 & 1\\
                             \end{array}\right ),
                         \left (\begin{array}{cc}
                               1 &    b\lambda     \\
                               0 & 1\\
                             \end{array}\right )\in K.\eqno(4.1)$$
 Since gcd$\,(a,b)=1$,
 there exists $m,n\in \Bbb Z$ such that $ am+bn=1$. As a consequence,
  we have
$$ \,\,\,\,  T =
 \left (\begin{array}{cc}
                               1 &    a\lambda      \\
                               0 & 1\\
                             \end{array}\right )^m
                         \left (\begin{array}{cc}
                               1 &    b\lambda     \\
                               0 & 1\\
                             \end{array}\right )^n =
\left (\begin{array}{cc}
                               1 &    \lambda     \\
                               0 & 1\\
                             \end{array}\right )
                             \in K.\eqno(4.2)$$
 Since $H_5$ is generated by $S$ and $T$,
 every element in $H_5$ is a word $w(S,T)$.
  Since  $T \in K$,  $SK = KS$ ($K$ is normal),
$ w(S,  T)K $ is  $S^iK$  for some $i$ $(S$ has order 4 in $H_5$). Hence
   the index
     of $K$ in $H_5$ is either 1,  2 or 4.
      This implies that $(ST)^4 \in K$.
      Since the order of $ST$ is 5, we conclude that
       $ST \in K$.
    It follows that
      $ST  ,   T
       \in K.$
      As a consequence,
       $ S\in K.$ Hence $K = H_5$.
        The lemma now can be proved by Second Isomorphism Theorem and the fact
         $H(\pi)\cap H(\tau) = H(\pi\tau)$.\qed

\smallskip
\noindent {\bf Remark 4.2.}  In the case $q=5$, $\Bbb Z[\lambda_5]$ is a principal ideal domain.
 As a consequence, every ideal of $\Bbb Z[\lambda_5]$ takes the form $(\tau)$ for some
  $\tau$.

\smallskip
\noindent {\bf Proposition 4.3.} {\em Let $p \in \Bbb N$ be an odd prime.
  Suppose that $H(p)$ contains an element $\sigma$
  such that {\tiny $\sigma \equiv
 \left (\begin{array}{cl}
                               1 &   p   \\
                               0 & 1\\
                             \end{array}\right )$} $(mod\, p^{2})$.
                              Then $[H(p^n) : H(p^{n+1})] = p^6$ for all $n \in \Bbb N$
                               and $[H_5 : H(p^n)] =
                      p^{6(n-1)}          [H_5 : H(p)]$.}

\smallskip

\noindent {\em Proof.} Let $X = \sigma^{p^{n-1}}$, $A= T^{p^n} \in H(p^n)$. It follows that
 $X \equiv $ {\tiny $
 \left (\begin{array}{cl}
                               1 &   p ^n  \\
                               0 & 1\\
                             \end{array}\right )$} $(mod\, p^{n+1})$
                              and
    $A=  $ {\tiny $
 \left (\begin{array}{cl}
                               1 &   p ^n \lambda  \\
                               0 & 1\\
                             \end{array}\right )$}. As a consequence, $H(p^n)$
 contains $A, SAS^{-1}, TSAS^{-1}T^{-1}$ and
$X, SXS^{-1}, TSXS^{-1}T^{-1}$. The above matrices modulo $p^{n+1}$ are given by
{\tiny $$\left (\begin{array}{cc}
   1 & p^n\lambda\\
   0 & 1
 \end{array}\right ),
 \left (\begin{array}{cc}
   1 & 0\\
  - p^n\lambda & 1
 \end{array}\right ),
 \left  (\begin{array}{cc}
   1-p^n\lambda^2  & p^n \lambda^3\\
   -p^n\lambda& 1+p^n\lambda^2
 \end{array}\right ),
 \left (\begin{array}{cc}
   1 & p^n\\
   0 & 1
 \end{array}\right ),
 \left (\begin{array}{cc}
   1 & 0\\
  -p^n & 1
 \end{array}\right ),
 \left  (\begin{array}{cc}
   1-p^n \lambda & p^n\lambda^2 \\
  - p^n & 1+p^n\lambda
 \end{array}\right ).$$}
Applying Lemma A in Appendix A, they generate an elementary abelian group of order
 $p^6$. Hence  $[H(p^n) : H(p^{n+1})] \ge  p^6$. By (2.2),
 $[H(p^n) : H(p^{n+1})] = p^6$.\qed

 \smallskip
\noindent {\bf Corollary 4.3.} {\em
  Suppose that $H(4)$ contains an element $\sigma$
  such that {\tiny $\sigma \equiv
 \left (\begin{array}{cl}
                               1 &     4\\
                               0 & 1\\
                             \end{array}\right )$} $(mod\, 8)$.
                              Then $[H(2^n) : H(2^{n+1})] = 2^6$ for $n \ge 2$
  and $[H_5 : H(2^n)] =
                      2^{6(n-2)}          [H_5 : H(4)]$.}

\smallskip
\noindent {\bf Remark.} In the case $p=2$, $H(2)$ possesses no $\sigma$
 such that $\sigma\equiv $ {\tiny $\left (\begin{array}{cl}
                               1 &     2 \\
                               0 & 1\\
                             \end{array}\right )$} $(mod\, 4)$.
  Consequently, Corollary 4.3 cannot be improved to $H(2)$.

\section {The index  $[H_5 :H(\pi)]$ where $N(\pi)$ is a
 power of a prime}

 Throughout the section, $N(\pi)$ is a power of a rational prime $p$.
  There are five  cases to consider, (i) $p=5$, (ii) $ p=2, $
  (iii) $ p = 3$,
   (iv) $p\ne 3$ is  of the form
   $10k\pm 3$, (v)  $ p$ is  of the form
   $10k\pm 1$.

 \subsection {$N(\pi)$ is a power of 5} Let $\tau = 2 +\lambda. $  The case $p=5$ is slightly different as
$5 = \tau^2 \lambda^{-2}$ ramified totally in $\Bbb Z[\lambda]$. As a consequence,
$ \pi = \tau^n$ for some $n $.
 We shall first determine the order of $H(5^n)/H( 5^{n+1})$.
  Lemma 3.3 ($p= 5$, $m = 5$)
 implies that  $S \in H(5)$, where
 $$
 \left (\begin{array}{cc}
   1 & 12\cdot 5 \\
   0 & 1
 \end{array}\right )
\equiv
 \left (\begin{array}{cc}
   1 & 20\lambda ^5\\
   0& 1
 \end{array}\right )\equiv S\,\,(mod\,\,25).\eqno(5.1)$$
  By Proposition 4.3,
  $[H(5^m) : H(5^{m+1})] = 5^6$ if $m \ge 1$.  By (2.2),
  $[H(\tau ^m) : H(\tau^{m+1})] \le 5^3$.  Hence
   $$[H(\tau ^m) : H(\tau^{m+1})] = 5^3 \mbox{ for } m\ge 2.\eqno(5.2)$$

\smallskip
    \noindent
    In the case $m=1$,
  applying  our results in subsection 3.3, $H( \lambda +2)$ contains the following matrix.
$$a =\left (\begin{array}{cc}
  -11\lambda - 6  & 10\lambda+5 \\
  4\lambda +3  & -4\lambda -2
\end{array}\right )= T ^{-2}
\left (\begin{array}{cc}
  3 \lambda +2 &-2\lambda-3\\
  4\lambda +3  & -4\lambda -2
\end{array}\right )\in H( \lambda +2).\eqno(5.3)$$
 Let $J=$ {\tiny $
 \left (\begin{array}{cc}
   0 & 1\\
   1& 0
 \end{array}\right )$}. Then $ J \in Aut H_5$.
  Let  $\Delta = \{ a, b=SaS^{-1},
 c= JaJ^{-1}\} \subseteq H(\tau)$. The matrices in
$\Delta$ modulo $5$ are  given by
$$
a\equiv I + (\lambda+2)\left (\begin{array}{cc}
  4 & 0 \\
  4  & 1
\end{array}\right ),
b\equiv I + (\lambda+2)\left (\begin{array}{cc}
  1 & 1 \\
  0  & 4
\end{array}\right ),
c\equiv  I + (\lambda+2)\left (\begin{array}{cc}
  1 & 4 \\
  0& 4
\end{array}\right ).\eqno(5.4)$$
    $\Delta$ generates a group of order $5^3$ modulo 5.
       By (2.2), $H(\tau)/H(5)  \cong \left  <\Delta\right >$ and
        $|H( \tau)/H(5)| =5^3$. In summary, one has
     $$[ H( \tau^{m}) : H( \tau^{m+1})] = 5^{3}
    \mbox{  for } m \ge 1.\eqno(5.5)$$
Recall that $[H_5 : H(\tau)] = 120 $ and that $N(x)$ is the absolute norm of $x$. It follows from (5.5) that
$$[H_5 : H(\tau ^m)] = 120\cdot 5^ {3(m-1)} = N(\tau^m)^3 (1-N(\tau)^{-2})
\mbox{ for } m\ge 1  .\eqno (5.6)$$

 \subsection {$N(\pi)$ is a power of 2} Since  $2$ is a prime in $\Bbb Z[\lambda]$,
  $\pi = 2^m$. Note first that $[H_5 : H(2)] = 10.$
 We shall now determine the index $[G(2):G(4)]$.
 Applying Lemma 3.1, a set of  generators of $G(2)$ is given by
{\small $$\Omega_2 =\left \{\left (\begin{array}{cc}
   1 & 2\lambda\\
   0 & 1
 \end{array}\right ),
 \left (\begin{array}{cc}
   1 & 0\\
   2\lambda & 1
 \end{array}\right ),
 \left (\begin{array}{cc}
   1+2\lambda & 2+ 2\lambda\\
   2\lambda & 1 +2\lambda
 \end{array}\right ),
 \left (\begin{array}{cc}
   1+2\lambda & 2\lambda\\
   2 +2\lambda & 1+2\lambda
 \end{array}\right )\right \}.\eqno(5.7)$$}

\noindent  Since $\Omega_2$ generates $G(2)$,
 one can show by direct calculation that
   $G(2)/G(4)$ is elementary abelian of order $2^4$.
   Since $-I \in H(2)-H(4)$,
   one has  $[H( 2) : H(4)] = 2^5$.
 We now study  $[H(2^n):H(2^{n+1})]$ ($n\ge 2$). By (2.2),
  $[H(2^n) :H(2^{n+1}]\le 2^6$. By our results in subsection 3.3,
 {\small   $$ 
    \left (\begin{array}{cc}
   1&  0\\
   -4\lambda & 1
 \end{array}\right )
    \left (\begin{array}{cc}
   1&  -4\lambda \\
   0 & 1
 \end{array}\right )
     \left (\begin{array}{cc}
   1+2\lambda & 2+ 2\lambda\\
   2\lambda & 1 +2\lambda
 \end{array}\right )^2
    \equiv  \left (\begin{array}{cc}
   1 & 4\\
   0 & 1
 \end{array}\right ) \in H(4)/H(8).\eqno(5.8)$$}

  \noindent By Corollary 4.4,
  $$
  [H_5 : H(2)]=10,
 [H_5 : H(2^n )]= 5\cdot 2^{6(n-1)}  \mbox{ for } n \ge 2.\eqno(5.9)$$

\subsection { $N(\pi)$ is a power of 3}
  Since   $3 \in \Bbb Z[\lambda]$  is a prime, $\pi = 3^m$.
  Note first that $[H_3 : H(3)] = 120$.
  Applying our results in subsection 3.3,
 $$ A= \left (\begin{array}{cc}
   \lambda & \lambda+2\\
   \lambda & 2\lambda+1
 \end{array}\right )
 \left (\begin{array}{cc}
   -1 & \lambda \\
  -2\lambda & 2\lambda+1
 \end{array}\right )^{-1}
 =
\left (\begin{array}{cc}
   9\lambda+4 & -2\lambda -3\\
  9\lambda +6& -3\lambda-2
 \end{array}\right ) \in H_5
 .\eqno(5.10)$$
  Direct calculation shows that
  $  E_3 = A^{-3} \equiv {\tiny
  \left (\begin{array}{cc}
   1 & 3\\
   0 & 1\\
 \end{array}\right )}\,\,(\mbox{mod } 9).$
  Applying  Proposition 4.3, we have
$$[H_5 : H(3^m)]= 120 \cdot 3^{6(m-1)}  \mbox{ for } m\ge 1 .\eqno(5.11)$$


\subsection {$N(\pi)$ is a power of $p$, where $3\ne p$ takes the form $10k\pm 3$}
 It follows that $p$ is a prime in $\Bbb Z[\lambda]$ and $\pi = p^n$.
Note first that $[H_5 : H(p)] = |SL(2, p^2)|.$
 By Lemma 3.3 ($m=p$),
     $$  T^{-20p}S \equiv
 \left (\begin{array}{cc}
   1 &  12p \\
   0 & 1
 \end{array}\right )
\in H( p)/H( p^2 ).\eqno(5.12)$$
By Proposition 4.3,
$$[H_5 : H(p^n)] = p^{6(n-1)} |SL(2, p^2)| = N(p^n)^3 (1-N(p)^{-2}).\eqno (5.13)$$

\subsection {$N(\pi)$ is a power of $p$, where $p$ takes the form $10k\pm 1$}
 $p$ is not a prime in $\Bbb Z[\lambda]$ and $ p = \tau \sigma$ where $\tau$ and $\sigma$
  are primes in $\Bbb Z[\lambda]$. Consequently,
   $\pi = p^r\tau^s$.
By Lemma 3.3 ($m=p$),
       $$T^{-20p}S\equiv
 \left (\begin{array}{cc}
   1 &  12p \\
   0 & 1
 \end{array}\right )
\in H(p)/H(p^2).\eqno(5.14)$$
By Proposition 4.3,
$[H( p^m) : H(5, p^{m+1})] = p^6$.
 It follows immediately from  (2.2)
   that $[H( \tau^mp^n) : H(5, \tau ^{m+1}p^{n})]
   = p^3 \mbox{ for } m\ge 1 . $
  In summary,
   $$
   [H( p^m) : H( p^{m+1})] = p^6,\,\,
   [H( \tau^mp^n) : H(\tau ^{m+1}p^{n})]
   = p^3  \mbox{ for } m\ge 1    . \eqno(5.15)$$

   \noindent {\bf Lemma 5.1.} {\em Let $p = 10k\pm 1 \in \Bbb N$ be a rational prime
   and let $p = \tau\sigma$
    where $\tau$ and $\sigma$ are primes in $\Bbb Z[\lambda]$.
   Then
   $[H_5 : H(p)] = [H_5 : H(\sigma)][H_5 : H(\tau)] = ((p-1)p(p+1))^2$.}

   \smallskip
   \noindent {\em Proof.} Let $H_0(p)$ be given as in (1.2).  Applying results of
     [CLLT] and [LLT2],
 \begin{enumerate}
 \item[(i)]  $[H_5 : H_0( p)] = (p+1)^2$  (Lemma 1 of [CLLT]),
\smallskip
 \item[(ii)]  $H(\sigma)\cap H(\tau)
  = H(p)$,  $[H_5 : H(\sigma)] =[H_5 : H(\tau)]= (p-1)p(p+1)$,
  \smallskip
 \item [(iii)]
  $G_5/G(\sigma ) \cong G_5/G(\tau) \cong PSL(2, p)$ is simple (Corollary 2 of [LLT2]).
  \end{enumerate}
  By (i) and (ii) of the above,
  $H(\sigma ) \ne H(\tau)$.
     By (iii) of the above, $H(\sigma)H(\tau)= H_5$.
    Applying  Second Isomorphism Theorem, one has
    $[H_5 : H(p)] = [H_5 : H(\sigma )][H_5 : H(\tau)].$
     Hence
   $[H_5 : H(p)] = ((p-1)p(p+1))^2  .$\qed

\smallskip
The index   $[H_5 : H(p^r\tau^s)]$
 can be determined by applying (5.15) and Lemma 5.1. In summary,    $$[H_5 : H( \pi)]=  N(\pi)^3
\prod_{P|(\pi)} (1-N(P)^{-2}\eqno(5.16)$$

\section {The Main Results}
Let $X \in \Bbb Z[\lambda]$.  Consider the decomposition $X =\prod x_i$,
 where $N(x_i)$ is the power of a rational prime $p_i$ such that
  gcd$\,(p_i, p_j) = 1$ for all $i \ne j$. By Lemma 4.1,
   $[H_5 : H(X)]
   = \prod [H_5 : H(x_i)]$. The index $[H_5 : H(x_i)]$  can be determined
    by applying (5.6), (5.9), (5.11), (5.13), and (5.16). In short, we have,
    if gcd$\,(N(\pi), 6)=1$, then
      $$[H_5 : H( 2^a3^b\pi)]= I_aJ_b\, N(\pi)^3
\prod_{P|(\pi)} (1-N(P)^{-2},\eqno(6.1)$$
where
the product is over the set of all prime ideals that divide $\pi$,
 $I_0 = 1,$ $I_1= 10$, $I_a = 5\cdot 2^{6(a-1)}$ if $a \ge 2$, $J_0=1$,
 $J_b = 120\cdot 3^{6(b-1)}$ if $b \ge 1$.
As a corollary of (6.1), one has the following.

\smallskip
\noindent {\bf Corollary 6.1.} {\em Let $A$ be an ideal. Then the natural homomorphism
  $f : H_5 \to  SL(2, \Bbb Z[\lambda ]/A)$ is surjective
   if and only if the norm of $A$ is prime to $6$. }

\section {Application : The commutator subgroup of $H_5$ is not congruence}
\subsection{}
Applying our results in subsection 5.1, we have the following.
\begin{enumerate}
\item[(i)] $H_5/H(5) \cong  H(\lambda +2)/H( 5 )  \cdot   H_5/H( \lambda+2)\cong  E_{5^3} SL(2,5)$, where $ H( \lambda +2)/H( 5 )\cong E_{3^5}\cong
 \Bbb Z_5 \times \Bbb Z_5 \times \Bbb Z_5$ is the elementary
 abelian group of order $5^3$ and $ H_5/H( \lambda+2) \cong SL(2, 5) $,
  $H(\lambda +2)/H(5) \cong \left < a, b, c\right > = \left < \Delta\right >$, where
$$
a\equiv I + (\lambda+2)\left (\begin{array}{cc}
  4 & 0 \\
  4  & 1
\end{array}\right ),
b\equiv I + (\lambda+2)\left (\begin{array}{cc}
  1 & 1 \\
  0  & 4
\end{array}\right ),
c\equiv  I + (\lambda+2)\left (\begin{array}{cc}
  1 & 4 \\
  0& 4
\end{array}\right ).\eqno(7.1)$$
 \end{enumerate}

\subsection{}$H_5^5$ and $H_5'$ are not congruence. Denote by $H_5^5$ the subgroup of
 $H_5$ generated by all the elements of the forms $x^5$, $x\in H_5$. $H_5^5$ is known as
  the power subgroup of $H_5$.
  It is clear that $H_5^5$ is normal and that $S, T^5 \in H_5^5$. Further, the index
   of $H_5^5$ in $H_5$ is 5.
  Since $H_5/H_5^5$ is abelian, $H_5^5$ contains the commutator
    subgroup $H_5'$.

\smallskip
\noindent {\bf Lemma 7.1.} {\em Suppose that $H_5^5$ is congruence. Then
 $H(5) \subseteq H_5^5$.}
 \smallskip

 \noindent {\em Proof.} Suppose that $H_5^5$ is congruence. Then $H(5^m r)\subseteq
  H_5^5$ for some $5^mr \in \Bbb N$, where $ m \ge 1$, gcd$\, (5, r)=1$.
  Let $K$ be the smallest normal subgroup that contains $H(5^m r)$ and $T^{5^m}$.
   It is clear that $K \subseteq H(5^m)\cap H_5^5$. Since
    gcd$\,(5, r)=1$, $T^{5^m}\in K $ and
     $T^r\in H(r)$, $$T \in KH(r) \subseteq H(5^m)H(r).\eqno (7.2)$$ Note that
     $KH(r)$ is a normal subgroup of $H_5$ and that $T\in KH(r)$. Applying the proof of Lemma 4.1,
      one has $H_5 = KH(r)$. Note that  $H(5^m) H(r) = H_5$ and that  $H(5^m)\cap H(r)
        = H(5^mr)$.
      By second isomorphism theorem,
  {\small     $$K/H(5^mr)\times H(r)/H(5^mr)=
      KH(r)/H(5^mr) = H_5/H(5^mr)  = H(5^m)/H(5^mr)\times H(r)/H(5^mr).$$}
      Hence  $K = H(5^m)$.
       As a consequence, $H(5^m) =K  \subseteq H_5^5$.
        Let $m$ be the smallest positive integer such that $H(5^m) \subseteq H_5^5$. Suppose
         that $m \ge 2$. Applying our main result (6.1),
          $$H_5/H(5^m) \cong SL(2, \Bbb Z[\lambda]/5^m) \cong SL(2, \Bbb Z[\lambda])/L(5^m)
          \eqno(7.3)$$
           (see Section 2 for notation). Hence
           $H_5/H(5^m)$ and $SL(2, \Bbb Z[\lambda])/L(5^m)$ have the
            same coset representatives. In particular, there exists $ \sigma \in H_5$
             such that
             $$\sigma H(5^m) =
             \left (\begin{array}{cc}
  1 & 1 \\
  0& 1
\end{array}\right ) H(5^m).\eqno(7.4)$$
Since $m \ge 2$, $\sigma^{5^{m-1}} \in H_5^5$. As a consequence, $H_5^5/H(5^m)$
 contains the following elements.
{\small $$
 T^{5^{m-1}}=
  \left (\begin{array}{cc}
  1 & 5^{m-1}\lambda \\
  0& 1
\end{array}\right ),
\sigma ^{5^{m-1}}\equiv
  \left (\begin{array}{cc}
  1 & 5^{m-1}\\
  0& 1
\end{array}\right )  \in H(5^{m-1})/H(5^m)\cap H_5^5 /H(5^m).\eqno (7.5)$$}
Now we consider the
 subgroup $V$  of $H_5^5/H(5^m)$ generated by $ T^{5^{m-1}}$ and $\sigma ^{5^{m-1}}$.
  Applying the proof of Proposition 4.3, $V = H(5^{m-1})/H(5^m)$.
   In particular, $H(5^{m-1}) \subseteq H_5^5$. This contradicts the
    minimality of $m$. Hence $m=1$ and $H(5) \subseteq H_5^5$.\qed

\smallskip
\noindent {\bf Proposition 7.2.} {\em $H_5^5$  and $H_5'$  are  not congruence.}

\medskip
\noindent {\em Proof.}  Since $H_5' \subseteq H_5^5$,
 it suffices
 to show that $H_5^5$ is not congruence.
Suppose that $H_5^5$ is congruence. By Lemma 7.1,
 $H(5) \subseteq H_5^5.$ Since $H_5/H( \lambda+2)\cong SL(2,5)$ has no normal
  subgroup of index 5
  and $H_5^5$ has index 5 in $H_5$, $H(\lambda+2)$ is not a subgroup of $H_5^5$.
    This implies that $H_5^5 H(\lambda+2) = H_5$.  By Second Isomorphism Theorem,
     $|H( \lambda+2)/[H_5^5 \cap H(\lambda+2)]| =5$
      and $|[H_5^5 \cap H( \lambda+2)]/H(5)|=5^2$. Note that
        $ E_{5^3} SL(2, 5) \cong H_5/H( 5)$ acts on
      $D = [H_5^5 \cap H(\lambda+2)]/H(5) \cong \Bbb Z_5 \times \Bbb Z_5$
       by conjugation.
       Note also that $D$ is a subgroup of $\left <\Delta \right > $ (see (7.1)).
    Recall that
    $$ J = \left (\begin{array}{cc}
  0 & 1 \\
  1& 0
\end{array}\right )\in \mbox{ Aut}\,H_5.\eqno(7.7)$$

\smallskip
\noindent
$D$ is invariant under the conjugation of
 $J$ and every element of $H_5$ (in particular, $S$ and $T$).
 However, one sees by direct calculation that the only nontrivial subgroup of $\left <\Delta \right >$
  invariant under $J$, $S$,    and $T$ is $\left <\Delta\right >$ itself (see Appendix B).
  A contradiction.
    Hence $H_5^5$ is not congruence. \qed

 \section {Appendix A}
\noindent
 {\bf Lemma A.} {\em  Let $p\in \Bbb N$ be a prime.
   Then $\Omega$ modulo $p^{n+1}$ generates a
   group of order $p^6$, where $\Omega $ is given as follows.}
 $$\Omega \equiv \left\{
  \left (
\begin{array}{cc}
1 &  p^n  \lambda ^i   \\
0 & 1 \\
\end{array}
\right ), \left (
\begin{array}{cc}
1 &  0\\
-p^n\lambda  ^i& 1 \\
\end{array}
\right ) ,\left (
\begin{array}{cc}
1 -p^n \lambda  ^{i+1}     &  p^n \lambda   ^{i+2}\\
 - p^n  \lambda  ^i& 1 + p^n  \lambda ^{i+1}\\
\end{array}
\right ) : i = 0, 1 \right \}.\eqno(A1)$$

\smallskip
\noindent {\em Proof.}
 Put the matrices in $(A1)$ into the forms $I +p^nU$ and  $I+ p^n V$. One sees easily that
 $$(I+p^nU)(I+p^nV)\equiv I +p^n(U+V)\,\,\, (\mbox{mod} \,\,p^{n+1}).\eqno(A2)$$
 Hence $\Omega$ modulo $p^{n+1}$ generates an abelian group. Note that
  $(B2)$  makes the multiplication of $I +p^nU$ and  $I+ p^n V$ into the addition of
   $U$ and $V$.
  In order to show $\Omega$ generates a group of order $p^6$ modulo $p^{n+1}$, we
   consider the following groups.
{\small $$ M =  \left < X_i =
  \left (
\begin{array}{cc}
1 &  p^n  \lambda ^i   \\
0 & 1 \\
\end{array}
\right ),  Y_i = \left (
\begin{array}{cc}
1 &  0\\
-p^n\lambda  ^i& 1 \\
\end{array}
\right )
 \right > ,\,\,
N = \left < Z_i =
\left (
\begin{array}{cc}
1 -p^n \lambda  ^{i+1}     &  p^n \lambda   ^{i+2}\\
 - p^n  \lambda  ^i& 1 + p^n  \lambda ^{i+1}\\
\end{array}
\right )\right > .$$}
It is easy to see that
 $M$ and $N$ are abelian groups of order $p^4$ and $p^2$ respectively.
Applying $(A2)$ and the fact that $N$ is abelian,  elements in $N$ take the following simple form
$$
Z_0^{c_0} Z_1^{c_1} \equiv
\left (
\begin{array}{lr}
1 -p^n\sum _{i=0}^{1} c_i\lambda^{i+1}    &   p^n\sum _{i=0}^{1} c_i\lambda^{i+2}                   \\
-p^n\sum _{i=0}^{1} c_i\lambda  ^{i}  & 1 + p^n\sum _{i=0}^{1} c_i\lambda^{i+1}  \\
\end{array}
\right ).\eqno(A3)$$
Note that we may assume that $0\le c_i\le p-1$.
 Similar to $(A3)$, elements in $M$ take the form
  $$  X_0^{a_0} X_1^{a_1}Y_0^{b_0}Y_1^{b_1}\equiv\left (
\begin{array}{cc}
1    &   p^n\sum _{i=0}^{1} a_i\lambda ^{i}                 \\
-p^n\sum _{i=0}^{1} b_i\lambda ^{i}  & 1 \\
\end{array}
\right ).\eqno(A4)$$
Suppose that
$ Z_0^{c_0} Z_1^{c_1} \equiv  X_0^{a_0} X_1^{a_1}Y_0^{b_0}Y_1^{b_1}$
 modulo $p^{n+1}$.
 An easy study of the (22)-entries of $(A3)$ and $(A4)$
 implies that
 $ 1 + p^n\sum _{i=0}^{1} c_i\lambda^{i+1}\equiv 1$ (mod $p^{n+1})$.
 Hence $c_0 =c_1=0$. As a consequence,
$M \cap N=\{1\}$. Hence
 $|\left < \Omega \right > | = |M||N| = p^6.$\qed

\medskip
\section {Appendix B}

\noindent
 {\bf Lemma B1.} {\em  Let $\pi = \lambda+2$ and let
$\Delta =\{a, b,c\}$, where
$a,b,c$ are given as in $(7.1)$. Then
the only nontrivial subgroup of
$\left <\Delta\right >$  invariant under the action of $S, T$ and $J$
 is $\left <\Delta\right >$ }.

 \smallskip
 \noindent {\em Proof.}
 Since $(I+\pi U)(I+\pi V) \equiv I+\pi(U+V)$ (mod 5), multiplication
 of $(I+\pi U)(I+\pi V)$ can be transformed into addition of $U$ and $V$.
  This makes the multiplication of
  matrices $a$, $b$, and $c$ easy. Consequently, one has
$${\tiny  r
=(ac)(ab) \equiv I + \pi
 \left (\begin{array}{cc}
   0 & 0\\
   3 & 0
 \end{array}\right ),
 s= (ac)(ab)^{-1} \equiv I + \pi
 \left (\begin{array}{cc}
   0 & 3\\
   0 & 0
 \end{array}\right ),
 t=bc  \equiv I + \pi
 \left (\begin{array}{cc}
   -3 & 0\\
   0 & 3
 \end{array}\right ).} $$
\noindent It is clear that $\left <\Delta \right >= \left <a,b,c\right >
=\left <r,s,t\right >$.
  Let $A, B \in G_5$. Set $A^B = BAB^{-1}$. Direct calculation shows that
  $$r^S = s^{-1}, r^T = rs^{-1}t^2, r^J= s,
  s^S=r^{-1}, s^T= s, s^J=r,
  t^S=t^{-1},t^T=st, t^J=t^{-1}.\eqno(B1)$$
   Denoted by $M$ a nontrivial subgroup  of $\left <r,s,t\right >$
  that is invariant under the conjugation of $J$, $S$ and $T$. Let $1\ne  \sigma = r^is^jt^k \in M$.
     One sees easily that
 \begin{enumerate}

\item[(i)] If $k \not \equiv 0$ (mod 5), without loss of generality,
 we may assume that $ k =1$.
  Then $\sigma^J \sigma^S= t^{-2} \in M$. It follows that $t\in M$.
   Hence $t^T =st \in M$. Consequently, $s\in M$. This implies
  $s^S = r^{-1} \in M$.
   In summary, $r, s, t\in M$.

 \item[(ii)]
   If $k \equiv 0$ (mod 5),
    then $\sigma$ takes the form $r^is^j$. Suppose that $i\equiv 0 $ (mod 5).
     Then $ 1\ne s^j\in M$. It follows that $s\in M$. Consequently, $r = s^T\in M$.
     Hence $rs^{-1}t^2 = r^T\in M$. As a consequence, $t \in M$.
      In summary, $r, s, t\in M$. In the case $i\not\equiv0$ (mod 5),
       we may assume that $i=1$. Hence $rs^j\in M$.
        It follows that $(rs^j)^T (rs^j)^{-1} = s^{-1}t^2 \in M$.
        Consequently, $(s^{-1}t^2)^T=  st^2 \in M$.
         This implies that  $(s^{-1}t^2)(  st^2)  = t^4 \in M$.
          Hence $t\in M$. One now sees easily that
           $r, s, t \in M$.

\end{enumerate}
\noindent Hence  the only nontrivial  subgroup of $\left< \Delta\right >$
 invariant under  $J$, $S$ and $T$ is $\left <\Delta\right>$.\qed

\section {Appendix C}
 Applying  results of [K] and [LLT1],
 $\{ S, ST\}$  (see (1.1) for notation) is a set of independent generators of $G_q$.
 It follows that
 the commutator subgroup $G_q'= [G_q, G_q]$ of $G_q$
  has index $2q$  and $G_q/G_q' \cong \Bbb Z_q \times \Bbb Z_2$.

  \smallskip

\noindent {\bf C1}. Let {\small $$ T=
 \left (\begin{array}{cc}
   1 & \sqrt 2\\
   0 & 1
 \end{array}\right ),
 L=
 \left (\begin{array}{cc}
   1 &  0\\
   -\sqrt 2 & 1
 \end{array}\right ),
  S=
 \left (\begin{array}{cc}
   0 & 1\\
   -1 & 0
 \end{array}\right ).\eqno(C1)$$}

\noindent  Then
 $\overline G =G_4/G(4) = \left < T, L, S\right >$ has order 32
  (Parson [P]). Let $\overline V$ be the subgroup of $G_4/G(4)$
  generated by $TL$ and $LT$. It follows by direct calculation that $\overline V$ is a
   normal subgroup of $G_4/G(4)$ and that $\overline G/\overline V \cong \Bbb Z_4
   \times \Bbb Z_2$. Let $V$ be the inverse image of $\overline V$ in $G_4$.
    Then $G_4/V  \cong   \Bbb Z_4
   \times \Bbb Z_2 $.  Since the commutator subgroup $[G_4 , G_4]$ of $G_4$ has index 8 in $G_4$,
    we conclude that $ [G_4 , G_4] = V$. Hence
      $[G_4 , G_4]$ is congruence as $G(4) \subseteq V = [G_4 , G_4]$.
 \medskip

\noindent {\bf C2}. Let
{\small $$ T=
 \left (\begin{array}{cc}
   1 & \sqrt 3\\
   0 & 1
 \end{array}\right ),
 L=
 \left (\begin{array}{cc}
   1 &  0\\
   -\sqrt 3 & 1
 \end{array}\right ),
  S=
 \left (\begin{array}{cc}
   0 & 1\\
   -1 & 0
 \end{array}\right ).\eqno(C2)$$}

\noindent    Applying results of Parson [P], $G_6/G(2) \cong D_{12}$ has order 12.
 Hence $G_6/G(2)$ has a normal subgroup $\overline U$ of index 2. Further,
 $$G_6/G(3) = \left < T, L,  S\,:\,
 T^3 = L^3=S^2 = 1, TL=LT, STS^{-1} = L
 \right >  \cong (\Bbb Z_3\times \Bbb Z_3)\cdot \Bbb Z_2.\eqno(C3)$$
 Let $\overline V$ be the subgroup of $\overline G= G_6/G(3)$ generated by $TL^{-1}$.
 Then $\overline V$ is normal in $\overline G$ and $\overline G/\overline V  \cong \Bbb Z_6$.
 Let $U$ and $V$ be the inverse images of $\overline U$ and $\overline V$ in $G_6$.
  By second isomorphism theorem, $G(6) \triangleleft  UV \triangleleft  G_6$ and
   $G/UV \cong \Bbb Z_6\times \Bbb Z_2$.
Since the commutator subgroup $[G_6 , G_6]$ of $G_6$ has index 12 in $G_6$,
    we conclude that $ [G_6 , G_6] = UV$. Hence
      $[G_6, G_6]$ is congruence as $G(6) \subseteq UV = [G_6 , G_6]$.

\smallskip
\noindent {\bf C3}. Let $\Gamma (n)$ be the principal congruence subgroup
 of the modular group $\Gamma = PSL(2, \Bbb Z)$ of level $n$.
  Then $\Gamma/\Gamma (2) \cong S_3$ has a normal subgroup of index 2 and
   $\Gamma/\Gamma (3) \cong A_4$ has a normal subgroup of index 3.
    It follows that $\overline \Gamma =\Gamma/\Gamma (6)$ has a normal subgroup
     $\overline V$ of index 6 and $\overline \Gamma /\overline V \cong \Bbb Z_6$.
 Let $V$ be the inverse image of $\overline V$ in $\Gamma $.
    Then $\Gamma /V  \cong   \Bbb Z_6$.  Since the commutator subgroup $\Gamma '$ of $\Gamma$ has index 6 in $\Gamma $,
    we conclude that $ \Gamma '= V$. Hence
      $\Gamma '$ is congruence as $\Gamma (6) \subseteq V = \Gamma '$.

\bigskip{\small

\noindent Cheng Lien Lang\\
\noindent Department of Mathematics, I-Shou  University, Kaohsiung, Taiwan,
Republic of China.

\noindent   \texttt{cllang@isu.edu.tw}

\smallskip
\noindent Mong Lung Lang \\
\noindent Singapore 669608, Republic of Singapore.

\noindent \texttt{lang2to46@gmail.com}}

\medskip


\end{document}